\documentstyle[amsfonts,12pt]{article}

\newcommand{\lon}{\longrightarrow}
\newcommand{\rar}{\rightarrow}

\newcommand{\p}{{\partial}}

\newcommand{\C}{{\Bbb C}}
\newcommand{\ot}{\otimes}

\newcommand{\Img}{\mbox{\rm Im}\,}
\newcommand{\Ker}{\mbox{\rm Ker}\,}
\newcommand{\Beq}{\begin{equation}}
\newcommand{\Eeq}{\end{equation}}
\newcommand{\Beqr}{\begin{eqnarray*}}
\newcommand{\Eeqr}{\end{eqnarray*}}
\newcommand{\Barr}{\begin{array}}
\newcommand{\Earr}{\begin{array}}
\newcommand{\f}{{\cal O}}

\newcommand{\cM}{{\cal M}}

\newcommand{\al}{\alpha}
\newcommand{\be}{\beta}
\newcommand{\ga}{\gamma}

\newcommand{\la}{\lambda}
\newcommand{\om}{\omega}

\newcommand{\sip}{\smallskip}
\newcommand{\bip}{\bigskip}

\begin{document}
\sloppy

 \title{Formality of canonical symplectic complexes \\ and Frobenius
manifolds}

\author{S.A.\ Merkulov}

\date{}
\maketitle
\sloppy

\begin{abstract}
It is shown that the de Rham complex of a symplectic manifold $M$
satisfying the hard Lefschetz condition is formal. Moreover, it is
shown that the differential Gerstenhaber-Batalin-Vilkoviski algebra associated
to such a symplectic structure gives rise, along the lines
explained in the papers of Barannikov and Kontsevich [alg-geom/9710032] and
Manin [math/9801006], to the structure of a Frobenius manifold on the de
Rham cohomology of $M$.
\end{abstract}

\bip

\begin{center}
{\bf \S 0. Introduction}
\end{center}
It was shown in \cite{Ba-Ko} (see also \cite{Ma} for detailed
exposition and proofs) that  the formal
moduli space of solutions to the Maurer-Cartan equations modulo gauge
equivalence associated to a very special class of differential
Gerstenhaber-Batalin-Vilkoviski (dGBV) algebras, carries a natural
structure of a Frobenius manifold.

\sip

To author's knowledge, only one example of such a special
dGBV algebra was known, the one
constructed out of the Dolbeault complex of an arbitrary Calabi-Yau
manifold by Barannikov and Kontsevich \cite{Ba-Ko}.

\sip

In this note we produce another example of a special dGBV-algebra, this
time the one associated with an arbitrary symplectic manifold $(M,\omega)$ satisfying the
hard Lefschetz condition which says that the cup product
$$
[\omega^k]: H^{m-k}(M) \lon H^{m+k}(M)
$$
is an isomorphism for any $k\leq m=\frac{1}{2}\dim M$.
Applying then the machinery developed in \cite{Ba-Ko,Ma} to the moduli
space of solutions of the associated Maurer-Cartan equation we get
a structure of a Frobenius manifold on the de Rham cohomology of $M$.

\bip


\begin{center}
{\bf \S 1. Formality of the de Rham complex}
\end{center}

Let $M$ be a
$2m$-dimensional manifold equipped with a symplectic 2-form $\omega$.
The associated $(2m|2m)$-dimensional supermanifold $\cM=\Pi
TM$, $\Pi$ being the parity change functor and $TM$ the tangent bundle to $M$,
comes  equipped canonically with an odd
vector field $d$ and a second order even differential operator $L^*:
\f_{\cM} \rar \f_{\cM}$, where $\f_{\cM}$ is the (complexified)
structure sheaf on
$\cM$. They are most easily described in a local coordinate chart
$(x^a, \psi^b=dx^b)$, $a,b=1,\ldots,2m$, on $\cM$,
$$
d= \sum_{a=1}^{2m} \psi^a \frac{\p}{\p x^a},
$$
and
$$
L^* =\sum_{a,b=1}^{2m} \omega^{ab} \frac{\p^2}{\p \psi^a\p \psi^b},
$$
where $\omega^{ab}$ is the $2m\times 2m$ matrix inverse to the matrix,
$\omega_{ab}$, of coefficients of $\omega$ in the basis $dx^a$.
Under the canonical isomorphism $\Gamma(\cM, \f_{\cM})=\Omega^*M$ the
vector field $d$ goes into the usual de Rham differential.

\bip

{\bf 1.1. Lemma} {\em The second order differential operator
$\Delta:= [L^*, d]$
satisfies $\Delta^2=0$ and $[\Delta,d]=0$.}

\sip

\noindent{\em Proof}. In a local coordinate chart,
$$
\Delta= \sum_{a,b}  \omega^{ab} \frac{\p^2}{\p \psi^a\p
x^b}
- \sum_{a,b,c} \frac{\p\omega^{ab}}{\p x^c} \psi^c
\frac{\p^2}{\p \psi^a \p \psi^b}.
$$
Under the assumption (without loss of generality) that
$x^a$ are Darboux coordinates the required
statements become obvious. $\Box$

\bip

The isomorphism $\Gamma(\cM, \f_{\cM})=\Omega^*M$ sends $\Delta$ into
a differential $\Delta: \Omega^* M \rar \Omega^*M$ of
degree -1 on forms.

\bip

{\bf 1.2. Remark.} Clearly, for any manifold $M$ and any section $\nu\in
\Gamma(M, \Lambda^2 TM)$ we can define operators $d$, $L^*$ and
$\Delta=[L^*, d]$ on $\Pi TM$ as above.
Koszul \cite{Kz} showed that  Lemma~1.1 still holds true if the
pair $(M,\nu)$ is a Poisson manifold.
He suggested to call the cohomology of the resulting complex $(\Omega^*M,
\Delta)$ the {\em canonical cohomology}. Brylinski \cite{Br} showed that
for a symplectic manifold the canonical and de Rham cohomologies coincide.
He also showed that $\Delta$, when viewed as a degree $-1$
differential on $\Omega^* M$,  satisfies
$$
\Delta|_{\Omega^k M} =(-1)^{k+1}*d*,
$$
where  $*: \Omega^k M \rar\Omega^{2m-k}$ is the symplectic analogue of the Hodge
duality operator
defined by the condition $\be\wedge (*\al)= \langle \be, \al \rangle
\omega^m/m!$, with $\langle\, , \, \rangle$ being  the pairing between
$k$-forms induced by the symplectic form. This star operator satisfies
$*(*\al)=\al$ and $\be\wedge (*\al)= (*\be)\wedge \al$.

\bip

{\bf 1.3. Symplectic harmonic forms.}  A differential form $\al\in
\Omega^* M$ is called {\em symplectic harmonic} if it satisfies
$d\al =\Delta \al =0$. Mathieu \cite{Math} proved that the following
three statements are equivalent
\begin{itemize}
\item[(i)] the symplectic manifold $M$ satisfies the Hard Lefschetz
condition;
\item[(ii)] the morphism of differential complexes $(\Omega^*M, \Delta)
\rar (\Omega^* M/d\Omega^* M, \Delta)$ induces an isomorphism in
cohomology;
\item[(iii)] any class in the de Rham cohomology $H^*(M,\C)$ contains a
symplectic harmonic representative.
\end{itemize}

\bip

We use these results to prove the following

\sip

{\bf 1.4. Proposition.} {\em Let $M$ be a symplectic manifold satisfying
the Hard Lefschetz condition. Then the differentials $d,\Delta: \Omega^*M
\rar \Omega^*M$ satisfy}
$$
\Img d\Delta= \Img d \cap \Ker \Delta = \Img \Delta \cap \Ker d.
$$

\sip

{\em Proof.} It follows immediately from 1.3(ii) that
$\Img d \cap \Img \Delta= \Img d \cap \Ker \Delta = \Img \Delta \cap \Ker
d$. Thus it remains to show that
$\Img d \cap \Img \Delta= \Img d\Delta$
which will follow from the following

\sip

{\sc Claim.} {\em For any $p$-form $\al_p$ such that
$\al_p=d\ga_{p-1}= \Delta \be_{p+1}$ for some $\ga_{p-1}\in
\Omega^{p-1}M$ and $\be_{p+1}\in \Omega^{p+1}M$
 there exists a $p$-form $\tau_p$ such
that $\al_p=d\Delta \tau_p$. }

\sip

We shall prove this Claim by induction. It is trivially true for $p=2m$
(and $p=0$).
Let us show that it is true for $p=2m-1$.
Since $d\be_{2m}$ is trivially $0$, then, by 1.3(iii), there is a
representation $\be_{2m}= \be_{2m}^0 + d\tau_{2m-1}$ for some
$\tau_{2m-1}\in \Omega^{2m-1}M$ and $\be_{2m}^0\in \Omega^{2m}M$
satisfying $\Delta \be_{2m}^0=0$. Hence
$\al_{2m-1}=d\Delta \tau_{2m-1}$.

\sip

Assume now that the Claim is true for $p=k+2$. Let us show that it is true
for $p=k$. If $\al_k=d\ga_{k-1}= \Delta \be_{k+1}$, then, due to the
fact that $d$ and $\Delta$ commute, $\al_{k+2}:=d\be_{k+1}\in \Ker \Delta$.
Since $\Img d \cap \Ker \Delta = \Img d \cap \Img \Delta$, $\al_{k+2}=d\beta_{k+1}=
\Delta \mu_{k+3}$ for some $\mu_{k+3}\in \Omega^{k+3} M$ and hence, by
the induction hypothesis,
$\al_{k+2}= d\Delta \nu_{k+2}$ for some $\nu_{k+2}\in \Omega^{k+2} M$.
Then $d(\be_{k+1} - \Delta \nu_{k+2})=0$ and, by 1.3(iii), there is a
decomposition
$$
\be_{k+1} = \be_{k+1}^0 + d\tau_{k} + \Delta \nu_{k+2}
$$
for some $k$-form $\tau_k$ and $(k+1)$-form $\be_{k+1}^0$ satisfying
$d\be_{k+1}^0=\Delta \be_{k+1}^0=0$. Thus $\al_k=\Delta \be_{k+1}=d\Delta \tau_k$.
This completes the proof of the Claim and hence of the Proposition.
$\Box$

\bip

A differential complex is called {\em formal}\, if it is
quasi-isomorphic to its cohomology.

\sip

{\bf 1.5. Theorem.} {\em The de Rham complex $(\Omega^* M, d)$ on a symplectic manifold
satisfying the Hard Lefschetz condition is formal}.

\sip

{\em Proof.} It follows immediately from Proposition 1.4 above and Lemma
5.4.1 in \cite{Ma} that the natural inclusion
$$
(\Ker \Delta, d) \lon (\Omega^*M, d)
$$
and the projection
$$
(\Ker \Delta, d) \lon (H^*(M,\C), 0)
$$
induced from the map $\Ker \Delta \rar
H^*(\Omega^*M, \Delta)=H^*(M,\C)$, are quasi-isomorphisms. $\Box$

\sip

\bip
\begin{center}
{\bf \S 2. dGBV algebra of a symplectic manifold}
\end{center}
In this section we plug in the data of \S 1 into the general machinery
developed in \cite{Ba-Ko} (see also \cite{Ma}) and  produce the structure of a Frobenius
manifold on the de Rham cohomology of a symplectic manifold satisfying
the Hard Lefschetz condition. We shall give only a very short outline of the
construction and refer to \cite{Ma} for full details.

\sip

Let $(M,\om)$ be a symplectic manifold. For a moment  we switch back
to the interpretation of $\Delta$ and $d$ as an odd second order
derivation and, respectively,  an odd vector field on the supermanifold
$\cM=\Pi TM$.

\sip

{\bf 2.1. Odd Poisson structure on $\cM$.} For any $f,g\in \f_{\cM}$ we define
the odd brackets
$$
[f \bullet g]= (-1)^{\tilde{f}} \Delta (fg) - (-1)^{\tilde{f}} \Delta (f)g
- a\Delta b.
$$
where\, $\tilde{}$\, stands for the parity of the kernel symbol.
It is not hard to check that the conditions $\tilde{\Delta}=1$ and $\Delta^2=0$
imply
\cite{Ma}
\begin{itemize}
\item[a)] odd anticommutativity: $[f\bullet g]
=-(-1)^{(\tilde{f}+1)(\tilde{g}+1)} [g\bullet f]$;
\item[b)]  odd Jacobi identity:
$$
[f\bullet [g\bullet h]] = [[f\bullet g]\bullet h] +
(-1)^{(\tilde{f}+1)(\tilde{g}+1)} [g\bullet [f, \bullet h]];
$$
\item[c)] odd Poisson identity: $[f\bullet gh] = [f\bullet g]h +
(-1)^{\tilde{g}(\tilde{f}+1)}g [f\bullet h]$;
\item[d)] two odd differentials:
\begin{eqnarray*}
\Delta [f\bullet g] & = & [\Delta f \bullet g] +
(-1)^{(\tilde{f}+1)}[f\bullet \Delta g], \\
d [f\bullet g] & = & [d f \bullet g] +
(-1)^{(\tilde{f}+1)}[f\bullet d g].
\end{eqnarray*}
\end{itemize}
Thus $(\Gamma(\cM,\f_{\cM})=\Omega^*M, \bullet, \Delta, d)$ is an odd Lie
superalgebra with two commuting differentials. Note, however, that the
roles of $d$ and $\Delta$ are  not symmetric: $d$ is a
derivation of the associative multiplicative structure in $\Omega^* M$,
while $\Delta$ is not. Such a structure is often called a {\em
differential Gerstenhaber-Batalin-Vilkoviski algebra}.

\bip

{\bf 2.2. A normalised solution to the Maurer-Cartan equation.}
From now on we assume that $M$ satisfies the Hard Lefschetz condition.
Let $[c_i]$ be a basis and $x^i$ the  associated linear coordinates
in $H^*(M,\C)$. We define
$K=\C[[x^i]]$ and consider the odd Lie superalgebra $(K\ot_{\C}
\Omega^*M, \bullet)$ equipped with the differentials $d_K=1 \ot
d$ and $\Delta_K= 1\ot \Delta$. It follows from Proposition~1.4 above
and Proposition 6.1.1 in \cite{Ma} that there exists a generic even
formal solution $\Gamma\in K\ot \Omega^* M$ to the Maurer-Cartan equation
$$
d\Gamma + \frac{1}{2}[\Gamma \bullet \Gamma]=0
$$
such that $\Gamma_0=0$, $\Gamma_1= \sum_{i} x^i c_i$ and $\Gamma_n\in
K\ot \Img \Delta$ for all $n\geq 2$, where $c_i$ is a symplectic harmonic
harmonic representative of $[c_i]$ (with $c_0=1$),
and $\Gamma_n$ is the homogeneous component of $\Gamma$ of degree $n$ in
$(x^i)$. Moreover, $\Gamma$ can be chosen in such a way that all $\Gamma_n$
for $n\geq 2$ do not depend on $x^0$.

\sip

The operator
$$
\begin{array}{rccc}
d_{\Gamma}: & K \ot \Omega^*M & \lon & K\ot \Omega^*M \\
& f & \lon & d_{\Gamma}f= d_Kf + [\Gamma \bullet f]
\end{array}
$$
commutes with $\Delta$ and satisfies $d_{\Gamma}^2=0$. Actually, all the
results of \S 1  hold true after the replacements $\Omega^* M \rar
K\ot \Omega^* M$, $\Delta\rar \Delta_K$ and $d \rar d_{\Gamma}$.

\bip

{\bf 2.3. Integral.} Since $M$ satisfies the Hard Lefschetz condition,
$H^{2m}(M,\C)= \C[\omega^m]$ and hence $M$ is compact.
Then the integral
$$
\begin{array}{rccc}
\int_M: & \Omega^* M & \lon & \C\\
& \la &\lon & \int_M \la := \la\cap [M]
\end{array}
$$
is well-defined.

\sip

{\bf 2.3.1. Lemma.} {\em For any $\al,\be \in \Omega^* M$},
\begin{eqnarray*}
\int_M d\al\wedge \be &=& (-1)^{\tilde{\al}+1}\int_M\al\wedge d\be, \\
\int_M \Delta\al\wedge \be &=& (-1)^{\tilde{\al}}\int_M\al\wedge \Delta\be. \\
\end{eqnarray*}

{\em Proof}. The first statement follows immediately from the Stokes
theorem, while the second one requires a small computation (in which we
assume, for definiteness, that $\al\in \Omega^k M$ and hence $\be\in
\Omega^{2m-k+1}M$):
\begin{eqnarray*}
\int_M \Delta\al \wedge \be &=&  (-1)^{k+1}\int_M (*\,d*\al) \wedge
\be\\
&=&  (-1)^{k+1}\int_M (d*\al) \wedge (*\be)\\
&=&  \int_M (*\al) \wedge (d*\be)\\
&=&  \int_M \al \wedge (*\,d*\be)\\
&=& (-1)^k \int_M \al \wedge \Delta\be. \ \ \  \Box\\
\end{eqnarray*}
\sip

{\bf 2.4. From symplectic structures to Frobenius manifolds.}
Consider a map
$$
\begin{array}{rccc}
\psi:& H_K:=K\ot H^*(M,\C) & \lon & K\ot \Omega^* M \\
& X & \lon & \bar{X}\Gamma
\end{array}
$$
which, by definition, acts on the basis vectors $[c_i]$ of $H_K$  as follows
$$
\psi([c_i])= \frac{\p \Gamma}{\p x^i}.
$$
Using the isomorphism $\Ker d_{\Gamma}/ \Img d_{\Gamma}= H_K$,
one introduces a supercommutative structure into $H_K$,
$$
\overline{X\circ Y}:= \overline{X}\Gamma\cdot \overline{Y}\Gamma \bmod \Img
d_{\Gamma}.
$$
From Proposition  1.4 and Lemma~2.3.1 it easily follows that the data
$(\Omega^* M, d, \Delta, \int_M)$ satisfies the Assumptions 1-3 of
Manin in \cite{Ma}. Then his Theorems 6.2.3, 6.4.1 and
Proposition~6.3.1 \cite{Ma} immediately imply that the above product
is potential,
$$
[c_i]\circ [c_j]= \sum_{k,l}\frac{\p^3 \Phi}{\p x^i \p x^j \p x^k} g^{kl}[c_l]
$$
with
$$
\Phi= \int_M \left(\frac{1}{6}\Gamma^3 - \frac{1}{2}dB \Delta B\right),
$$
associative and admits an Euler vector field. Here
$\Gamma= \Gamma_1 + \Delta B$ with $B_0=B_1=0$, and $g_{ij}= \int_M
[c_i]\wedge [c_j]$ is the standard Poincare metric.

\sip

Thus $H^*(M,\C)$
carries the structure of a Frobenius manifold.

\sip

{\em Acknowledgements.} I am very grateful to the Max Planck Institute
for Mathematics for hospitality and excellent working conditions. I am
especially obliged to Yu.I.\ Manin who introduced me into the area and
made a number of insightful comments on the preliminary version of
the paper. Helpful remarks of J.-L.\
Brylinski and  A.N.\ Tyurin are gratefully acknowledged.

\bip

\noindent\mbox{\small Department of Mathematics, Glasgow University}

\noindent\mbox{\small 15 University Gardens, Glasgow G12 8QW, UK}

\noindent\mbox{\small e-mail: sm@maths.gla.ac.uk}

\end{document}